\PassOptionsToPackage{dvipsnames}{xcolor}

\documentclass[a4paper, 11pt]{elsarticle}

\usepackage[utf8]{inputenc}

\usepackage[lmargin=1.2in, rmargin=1.2in, bmargin=1.5in, top=1.2in]{geometry}

\usepackage{hyperref}

\usepackage{amsmath, amsfonts, mathrsfs}

\usepackage{graphicx, float}
\usepackage{xcolor}
\usepackage{lipsum}

\def\bb #1{\mathbb{#1}}
\def\cal #1{\mathcal{#1}}

\DeclareMathOperator\erf{erf}
\DeclareMathOperator\erfc{erfc}

\newcommand{\Fo}{\textit{Fo}}
\newcommand{\Ste}{\textit{Ste}}

\numberwithin{equation}{section}
\def\figtype{png}

\title{Physics Informed Neural Networks for heat conduction with phase change}
\date{\today}

\author[lmrs]{Bahae-Eddine Madir}
\ead{bahae-eddine.madir@univ-rouen.fr}

\author[lmrs]{Francky Luddens\corref{ca}}
\ead{francky.luddens@univ-rouen.fr}

\author[larema]{Corentin  Lothod{\'e}}
\ead{corentin.lothode@univ-angers.fr}

\author[lmrs]{Ionut Danaila}
\ead{ionut.danaila@univ-rouen.fr}

\address[lmrs]{Univ Rouen Normandie, CNRS, Normandie Univ, LMRS UMR 6085, F-76000 Rouen, France}
\address[larema]{Univ Angers, CNRS, LAREMA - UMR 6093, SFR MATHSTIC, F-49000 Angers, France}

\cortext[ca]{Corresponding author.}

\journal{Computer Method in Applied Mechanics and Engineering}

\begin{document}

\begin{abstract}
We study numerical algorithms to solve a specific Partial Differential Equation (PDE), namely the Stefan problem, using Physics Informed Neural Networks (PINNs). This problem describes the heat propagation in a liquid-solid phase change system. It implies a heat equation and a discontinuity at the interface where the phase change occurs. In the context of PINNs, this model leads to difficulties in the learning process, especially near the interface of phase change. We present different strategies that can be used in this context. We illustrate our results and compare with classical solvers for PDEs (finite differences).
\end{abstract}


\maketitle

\section{Introduction}

When modeling heat transfer problems, depending on the context, there are various formulations. For a solid material involving a liquid-solid phase change, the problem reduces to the Stefan problem, taking into account only the conduction heat transfer mechanism (natural convection is ignored). This problem was studied extensively by the mathematical community, resulting in many theoretical results (see \cite{rubinstein1979stefan} and references therein). In the context of numerical simulation, several frameworks were developed to solve this problem (see the review \cite{voller1990fixed}).\medskip

In the context of Finite Element Method (FEM), different approaches were suggested, from enthalpy based methods (e.g. \cite{elliott1987error}) to front-tracking methods (e.g. \cite{li1983finite}). For other numerical frameworks, similar approaches have been developed, see for example \cite{muhieddine2009various} for Finite Volumes Method (FVM) and \cite{savovic2003finite} for Finite Differences (FD).\medskip

Recently, machine learning methods have experienced significant growth, largely due to the improvement of computing hardware, such as GPUs. Thanks to these advancements, machine learning methods have found applications in various scientific domains, including heat transfer problems \cite{jambunathan1996evaluating, liu2018data, kwon2020machine, tamaddon2020data}. The majority of these applications relies on the supervised learning approach, for which large databases obtained from experiments or numerical simulations are required to train a model. One of the main challenges is the potential insufficiency of the data acquired through experiments. Additionally, high-fidelity numerical simulators are inherently slow, restricting their applicability in generating extensive numerical data-sets. Conversely, low-fidelity numerical simulators are generally faster, allowing for the creation of larger data-sets, albeit with a trade-off of reduced model accuracy. To tackle this problem, the physics-informed neural networks (PINNs) framework was developed, and subsequently it was applied to a various fluid mechanics problems \cite{RAISSI2019686, raissi2019deep, jin2021nsfnets, mao2020physics}, as well as heat transfer problems \cite{wang2021reconstruction, lucor2021physics, cai2021flow, wang2021deep, 10.1115/1.4050542}. Instead of solving partial differential equations (PDEs) using classical simulation tools to generate training data-sets, the governing PDEs can be used directly to train models.  One key advantage of such an approach is the efficiency in data assimilation problems where the PDE, initial or boundary conditions are not well known, but a finite number of sparse measurements inside the domain of interest is available. PINNs are also efficient when dealing with a PDE with parametric setting or high-dimensional PDE \cite{grohs2018proof, raissi2024forward}.\medskip

In this article, we study the applicability of PINNs, to solve the Stefan problem. The article is organized as follows. In Sec. 2, we describe Stefan's problem and the enthalpy formulation; we solve the problem using a finite difference method to generate a reference solution. In Sec. 3, we present a brief overview of the PINNs method and we use it to solve the problem for two different configurations. We study also the loss weighting effect on the model training. In Sec. 4, we present different strategies that can be used to improve the model accuracy. We conclude in Sec. 5 with a brief summary and outlook.

\section{The one-dimensional Stefan problem}

\label{sec:problemmethodo}

\subsection{Heating of a semi-infinite material subject to phase change}

The Stefan problem is defined on a semi-infinite 1d domain $x\ge 0$. The domain is filled with a pure material, subject to liquid-solid phase change. The dimensionless temperature and enthalpy are defined as:
$$\theta \leftarrow  \frac{T - T_f}{\delta T}, \qquad H \leftarrow \frac{H}{\rho\, c\, \delta T},$$ 
where the temperature scale for melting and solidification is $\delta T = \max(T_h - T_f, T_f - T_c)$, with $T_h$, $T_c$, and $T_f$ the hot, cold, and fusion temperatures, respectively. The constants $c$ and $\rho$ are the specific heat and the density (assumed to be the same in the liquid and the solid). The left wall, denoted by $\Gamma_l$ is kept at a constant temperature $\theta_l$, and the initial temperature at $t=0$ is $\theta_r$. The problem can be modeled as: 

\begin{equation}\label{eq:exact}
\begin{cases}
&\partial_t H - \Fo\,\partial_{xx}\theta = 0, \quad t>0,\ x>0,\\[1ex]
&\theta(t,0) = \theta_l, \\
&\lim_{x\to+\infty}\theta(t,x) = \theta_r,
\end{cases}
\end{equation}

where $\Fo = \alpha\, t_{\rm ref}/x^2_{\rm ref}$ is the Fourier number, with $\alpha$, $t_{\rm ref}$ and $x_{\rm ref}$ denoting the thermal diffusivity of the material, the reference time and reference length, respectively. The dimensionless enthalpy $H$ is defined as
\begin{equation}\label{eq:def_enthalpy}
H(t,x) = \begin{cases} \theta(t,x) &\textnormal{ if }\theta(t,x)\leq 0,\\
\theta(t,x) + \dfrac{1}{\Ste} &\textnormal{ if }\theta(t,x)> 0.\\[1ex]
\end{cases}
\end{equation}
The Stefan number is defined as $Ste = c\,\delta T/L$, with $L$ denoting the latent heat.\medskip

In this setting, the liquid-solid interface is a point $\lambda(t)$, and its motion is prescribed by the Stefan condition:
$$
\partial_x\theta_\text{right}(t,\lambda(t)) - \partial_x\theta_\text{left}(t,\lambda(t)) = \frac{1}{\Fo\,\Ste}\lambda'(t),
$$
where $\theta_\text{left}$ (resp. $\theta_\text{right}$) corresponds to $\theta_\text{ex}(t,x)$ for $x<\lambda(t)$ (resp. $x>\lambda(t)$).\medskip 

This condition is embedded in the enthalpy formulation~\eqref{eq:exact}. Using the ansatz $\lambda(t) = 2\lambda_0\sqrt{t\Fo}$, for some constant $\lambda_0$, the exact solution of the problem is given by:\medskip
\begin{equation}\label{eq:exact_solution}
\theta_\text{ex}(t, x) = \left\{
\begin{array}{lr}
\theta_l \left(1 - \erf\left(\lambda_0\right)^{-1}\erf\left(\frac{x}{2\sqrt{t\Fo}}\right)\right), & \text{ if } x \leq \lambda(t), \\[3ex]
\theta_r \left(1 - \erfc\left(\lambda_0\right)^{-1}\erfc\left(\frac{x}{2\sqrt{t\Fo}}\right)\right), &\text{ otherwise.}\\[2ex]
\end{array}\right.
\end{equation}
The Stefan condition is satisfied provided $\lambda_0$ is the solution of the nonlinear equation: 
$$
\lambda_0 - \frac{\Ste}{\sqrt{\pi}}e^{-\lambda_0^2}\left(\frac{\theta_r}{\erfc{(\lambda_0)}} + \frac{\theta_l}{\erf{(\lambda_0)}}\right) = 0.
$$

This solution is referred to as the \textit{exact solution} and will be used to set a suitable boundary condition on a finite domain.

\subsection{Restriction to a finite domain}

Numerically, we do not compute the solution on a semi-infinite domain. If the right boundary is kept at the temperature $\theta_{r}$, then the exact solution holds only for small times, for which the interface is far away from the boundary. In order to avoid using a large computational domain, we need to solve numerically the problem on a smaller domain, with a time-dependent Dirichlet condition that matches $\theta_\text{ex}$.\medskip

We consider a unit domain $D=[0, 1]$. The left and right walls, denoted by $\Gamma_l$ and $\Gamma_r$ are kept at a temperature $\theta_l$ and $\theta_\text{ex}(t,1)$ respectively. The evolution of the temperature is described by the following PDE:
\begin{align}
&\partial_t H(t, x) - \Fo\,\partial_{xx} \theta(t, x) = 0,\quad t>0,\quad x\in D, \label{eq:enthalpy_H}\\
&\theta(t, 0) = \theta_l , \\
&\theta(t, 1) = \theta_\text{ex}(t,1).
\end{align}

Noting that \eqref{eq:def_enthalpy} can be reformulated as $H = \theta + (1/\Ste)\varphi(\theta)$, where $\varphi$ is the Heaviside function representing the liquid fraction, the previous equation can be rewritten as
\begin{equation} \label{eq:enthalpy_T}
\partial_t \theta(t,x) - \Fo\,\partial_{xx} \theta(t,x) = -\frac{1}{\Ste}\partial_t \varphi(\theta).
\end{equation}

This equation is solved numerically in \cite{sadaka2020parallel, RAKOTONDRANDISA2020107188} using a regularized enthalpy. Let $\delta>0$ be a regularization parameter, and define
\begin{equation}\label{eq:phi_eps}
\varphi_\delta(\theta) = \frac12\left( 1 + \tanh\left(\frac\theta\delta\right)\right).
\end{equation}

Using the regularization $\varphi_\delta$ instead of $\varphi$, equation \eqref{eq:enthalpy_T} reads
\begin{equation}\label{eq:CN}
\partial_t\theta = r(\theta)\, \partial_{xx}\theta\qquad t \in [0.05, 1],\quad x \in [0, 1],
\end{equation}
where 
$$
r(\theta) = \Fo \left(1 + \frac{1}{\Ste}\varphi_\delta'(\theta) \right)^{-1}.
$$

\subsection{Reference solution}

\label{subsec:ref_sol}

As a regularization term has been used, $\theta_\text{ex}$ is no longer the solution of (2.9). In order to assess the accuracy of the neural network approach, we compute a reference solution of (2.9) using a finite differences scheme. Applying the Crank–Nicolson scheme~\cite{crank_nicolson_1947}  for time integration of \eqref{eq:CN} and centered finite differences for the space discretization, we obtain
\begin{equation}
\frac{\theta^{n+1}_i - \theta^{n}_i}{\Delta t} = \frac{1}{2}\left[r(\theta^{n+1}_i) \frac{\theta^{n+1}_{i + 1} - 2\theta^{n+1}_i + \theta^{n+1}_{i-1}}{(\Delta x)^2} + r(\theta^{n}_i) \frac{\theta^{n}_{i + 1} - 2\theta^{n}_i + \theta^{n}_{i-1}}{(\Delta x)^2} \right].
\end{equation}

\begin{figure}[H]
\centering
\includegraphics[scale=.65]{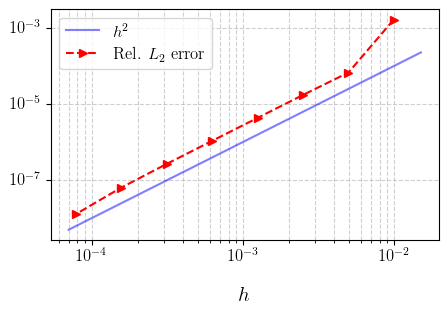}
\caption{Convergence of the finite differences scheme used to solve \eqref{eq:CN}. The $L_2$ error is calculated relative to a reference solution obtained with a step size $h_{\rm min} = 3.9\times 10^{-5}$, i.e. Rel. $L_2$ error: $\|\theta_h - \theta_{h_{\rm min}}\|_2/\|\theta_{h_{\rm min}}\|_2$.}\label{fig:1}
\end{figure}

Assuming that the solution $\theta^{n}$ is known, then $\theta^{n + 1}$ is  the solution of a non linear system of equations, which can be solved using the Newton–Raphson method \cite{BENISRAEL1966243}. Figure \ref{fig:1} shows the numerical convergence of the considered approach. By setting equal step size for time and space $h = \Delta t = \Delta x$, we obtain convergence of order $\cal O(h^2)$.\medskip

\section{Stefan problem with Physics-Informed Neural Networks}

\subsection{Physics-Informed Neural Networks}
In this section, we provide a brief overview of the Physics-Informed Neural Networks (PINNs) \cite{RAISSI2019686}, which is a method to approximate the solution to differential equations using neural networks.\medskip

A neural network with $L$ layers (or $L-1$ hidden layers), can be represented as the composition of $L$ parameterized functions $\ell^k (x, \Theta^k)$, where $x$ and $\Theta^k$ are the input and a set of parameters for the $k-$th layer, respectively. A feed-forward neural network in particular, apply for each layer $\ell^k$ a linear and nonlinear transformations to the inputs i.e.
$$
\forall k,\qquad \ell^k\left(x, \mathbf{W}^k; \mathbf{b}^k\right) = \Sigma^k\left(\mathbf{W}^k x + \mathbf{b}^k\right),
$$
where $\mathbf{W}^k\in \bb R^{N_{k-1}\times N_{k}}$ and $\mathbf{b}^k\in \bb R^{N_k}$ are the weights and the biases which constitute the neural network parameters. $N_k$ denotes the number of neurons in layer $k$ and the $\Sigma^k$ are nonlinear functions applied element-wisely, typically 
$$
\Sigma^k(x) = \left\{\begin{array}{ll}
\tanh(x), &\text{if}\quad  1\leq k<L,\\
x, &\text{if}\quad k = L.
\end{array}\right.
$$
In the context of finding the solution $u$ of a partial differential equation, a feedforward neural network can be considered as an approximate solution.\medskip

Consider PDEs taking the general form
\begin{equation}\label{eq:3.1}
\partial_t u + \cal N[u] = 0,\quad t\in [0, T],\quad x\in \Omega,
\end{equation}
subject to initial and boundary conditions
\begin{align}
&u(0, x) = g(x),\quad x\in \Omega,\label{eq:3.2}\\
&\cal B[u](t, x) = h(t, x),\quad t\in [0, T],\quad x\in \partial\Omega,\label{eq:3.3}
\end{align}
where $\cal N[\,\cdot\,]$ is a linear or nonlinear differential operator, and $\cal B[\,\cdot\,]$ is a boundary operator corresponding to Dirichlet, Neumann, Robin, or periodic boundary conditions.\medskip

By replacing the unknown solution $u(t, x)$ of \eqref{eq:3.1}, \eqref{eq:3.2} and \eqref{eq:3.3} with a neural network $u_\Theta (t, x)$, where $\Theta$ denotes all tunable parameters (weights and biases) of the network, the goal is to minimize the following composite loss function
\begin{equation}
\cal L(\Theta) = \omega_{0}\, \cal L_0(\Theta) + \omega_{b}\, \cal L_b(\Theta) + \omega_{r}\, \cal L_r(\Theta),
\end{equation}
where $\cal L_r$, $\cal L_0$ and $\cal L_b$ are loss (residual) terms corresponding to the PDE term \eqref{eq:3.1}, the initial condition \eqref{eq:3.2}, and the boundary condition \eqref{eq:3.3}, respectively:
\begin{align}
\cal L_r(\Theta) &= \frac{1}{N_r} \sum_{k=1}^{N_r} \left|\frac{\partial u_\Theta}{\partial t}(t_r^k, x_r^k) + \cal N[u_\Theta](t_r^k, x_r^k)\right|^2,\\
\cal L_0(\Theta) &= \frac{1}{N_0} \sum_{k=1}^{N_0} \left|u_\Theta(0, x_i^k) - g(x_i^k)\right|^2,\\
\cal L_b(\Theta) &= \frac{1}{N_b} \sum_{k=1}^{N_b} \left|\cal B[u_\Theta](t_b^k, x_b^k) - h(t_b^k, x_b^k)\right|^2.
\end{align}
The sets $\{x_i^k, g(x_i^k)\}_{k=1}^{N_0}$, $\{t_b^k, x_b^k, h(t_b^k, x_b^k)\}_{k=1}^{N_b}$, and $\{t_r^k, x_r^k\}_{k=1}^{N_r}$ represent $N_0$ initial training points, $N_b$ boundary training points, and $N_r$ collocation points, respectively. They are randomly selected using low discrepancy sequence techniques (also known as a quasi-Monte Carlo sample) such as Latin Hypercube \cite{stein1987large}, Sobol sequence \cite{sobol1967distribution} or Halton \cite{halton1960efficiency} sequence.\medskip

It is worth noting that loss terms include gradients with respect to input variables. It is thus possible to compute these exactly using automatic differentiation \cite{baydin2018automatic} at any point in the domain, without the need for manual computation. Moreover, hyper-parameters $\omega_0$, $\omega_b$ and $\omega_r$ enable the flexibility of assigning a different learning rate to each individual loss term in order to balance their interplay during model training. These weights may be user-specified or tuned automatically during training \cite{wang2021understanding, wang2022and}.

\subsection{Validation of the PINNs approach}

We consider the equation \eqref{eq:enthalpy_T} with $\Fo = 0.01$, $\Ste = 0.5$, a regularization $\delta = 0.05$ and a hot temperature $\theta_l = 1$. The initial condition is chosen to be the exact solution at the initial time $t=0.05$ (see \eqref{eq:exact_solution}). The specific system is summarized as follows:
\begin{subequations}
\label{eq:3.8}
\begin{align}
&\partial_t\theta - 0.01\, \partial_{xx}\theta + 2\, \partial_t\left[\varphi_\delta(\theta)\right] = 0, \qquad t\in [0.05, 1], \quad x\in [0, 1]\label{eq:3.8a}\\
&\theta(0.05, x) = \theta_\text{ex}(0.05, x),\label{eq:3.8b}\\
&\theta(t, 0) = \theta_l,\label{eq:3.8c}\\ 
&\theta(t, 1) = \theta_\text{ex}(t,1).\label{eq:3.8d}
\end{align}
\end{subequations}

Following the original work of Raissi et al. \cite{RAISSI2019686}, we represent the latent variable $\theta$ by a feed-forward neural network $\widehat{\theta}$ with 6 hidden layers and 20 neurons per hidden layer, we use the $\tanh$ as an activation function. In order to compute the loss functions, the Latin Hypercube Sampling method was adopted to generate $N_0 = 1024$ samples for the initial condition, $N_b = 256$ samples for the boundary conditions and $N_r = 10000$ collocation points in the spatio-temporal domain for the PDE term.\medskip

The model parameters are initialized using Xavier initialization \cite{glorot2010understanding}, and the training procedure of  the resulting PINN is done with full-batch gradient descent using the Adam optimizer \cite{kingma2014adam} for $100k$ iterations with an exponential decay learning rate $t\mapsto \eta\,\gamma^{t/\kappa}$ where $\eta = 10^{-3}$, $\gamma = 0.9$, $\kappa = 8000$ and $t$ is the training iteration. Table~\ref{tab:1} summarizes the hyperparameters of both the physics-informed neural network and the training.\bigskip

\begin{table}[H]
\center
\resizebox{13cm}{!}{
\begin{tabular}{l  p{1cm} l}
\hline
\multicolumn{3}{c}{Hyperparameters}\\%
\hline\rule{0mm}{8mm}\noindent
Architecture && 
\begin{minipage}{12cm}
Two inputs $(t, x)$, one output $\widehat{\theta}(t, x)$, six hidden layers of $20$ neurons with the activation function $\tanh$\medskip
\end{minipage}\\
Initialization        && Xavier\\
Optimizer            && Adam for $100k$ epoch\\
Learning rate   && Exponential decay learning rate $t\mapsto \eta\,\gamma^{t/\kappa}$, $\eta = 10^{-3}$, $\gamma = 0.9$, $\kappa = 8000$\\
Training data && $N_0 = 1024$, $N_b = 256$, $N_r = 10000$ using LHC sampling\\[.5ex]
\hline
\end{tabular}
}
\caption{Hyperparameters used to solve Eq. \eqref{eq:3.8}.}\label{tab:1}
\end{table}

To evaluate the accuracy of the physics-informed neural network, we compare the model's prediction $\widehat\theta$ with the reference solution $\theta$ (given in sec.~\ref{subsec:ref_sol}) at each time step, calculating the relative $L_2$ error $\|\widehat{\theta} - \theta\|_2/\|\theta\|_2$ on a grid of $500 \times 500$ points in the spatio-temporal domain. The prediction $\theta$ is derived from averaging multiple simulations initialized independently, allowing us to consider the effect of initialization. Sensitivity to initialization will be demonstrated by the standard deviation of the simulation ensemble.\medskip

\begin{figure}[H]
\centering
\includegraphics[scale=0.47]{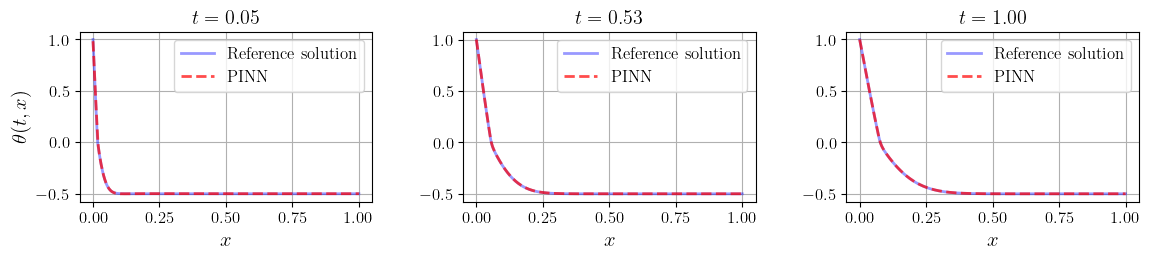}
\caption{Solution of~\eqref{eq:3.8}  at the end of training for $t=0.05$, $t=0.53$ and $t=1$.}\label{fig:13}
\end{figure}

In this setting, we test the baseline PINNs approach introduced in \cite{RAISSI2019686}, which sets $\omega_0 = \omega_b = \omega_r = 1$. From Fig. \ref{fig:13}, we notice that the neural network captured accurately the solution of the problem. The relative $L_2$ error at the end of training in this case is $3.175\times 10^{-3} \pm 0.0016$.

\begin{figure}[H]
\centering
\includegraphics[scale=.55]{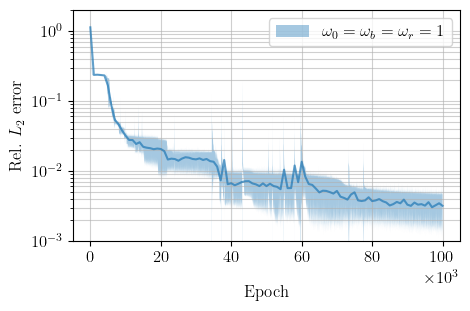}
\hspace{1ex}
\includegraphics[scale=.55]{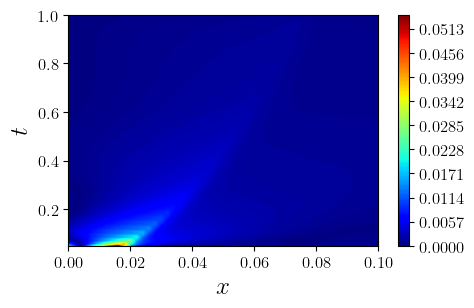}
\caption{Problem~\eqref{eq:3.8}. Left: relative $L_2$ error $\|\widehat\theta - \theta\|_2/\|\theta\|_2$ during training. Right: absolute error $|\widehat\theta - \theta|$ at the end of training (i.e. at epoch $=10^5$) .}\label{fig:12}
\end{figure}

\subsection{Influence of the Stefan number}

In this section, we assess the ability of PINNs to capture the phase change with smaller values of parameter $\Ste$. To this end, we change the Stefan number to the value $\Ste = 0.005$. Then Eq.~\eqref{eq:3.8a} becomes:
\begin{equation}\label{eq:3.9}
\partial_t\theta - 0.01\, \partial_{xx}\theta + 200\, \partial_t\left[\varphi_\delta(\theta)\right] = 0, \qquad t\in [0.05, 1], \quad x\in [0, 1].
\end{equation}

In this scenario, the disparity between the coefficients of the time derivative $\partial_t \theta$ and the spatial derivative $\partial_{xx}\theta$ is greater than it was previously. As mentioned in \cite{krishnapriyan2021characterizing}, this kind of discrepancy could potentially pose challenges during training, as the complexity of the loss landscape increases, making the minimization of the loss function more challenging.

For the model training, we maintain the same hyperparameters as previously. However, we introduce a supplementary adjustment by modifying the loss weights $\omega_0$, $\omega_b$, and $\omega_r$. We perform three test cases: we use $\omega_0 = \omega_b = \omega_r = 1$ as previously. Secondly, we set $\omega_0 = 100$ and $\omega_b = \omega_r = 1$ intended to enforce the neural network to satisfy the initial condition. Finally, we use the dynamical weighting algorithm suggested in \cite{wang2021understanding} to dynamically scale the weights of the loss function. The reweighting in this algorithm is performed so that the gradients of the loss terms of the initial and boundary conditions are approximately in the same range of values as the gradient of the PDE loss term. At each training step, e.g. iteration $(t+1)$, the estimates of $\omega_i$ and $\omega_b$ can be computed by:
\begin{equation}
\omega^{(t+1)} = (1 - \alpha)\,\omega^{(t)} + \alpha\, \frac{\max_\theta\{ |\nabla_\theta \cal L_r|\}}{\omega^{(t)}\overline{|\nabla_\theta \cal L|}},
\end{equation}
where $\max_\theta\{ |\nabla_\theta \cal L_r|\}$ is the maximum value attained by $|\nabla_\theta \cal L_r|$ and $\overline{|\nabla_\theta \cal L|}$ denotes the mean of $|\nabla_\theta \cal L|$. The initial values of the dynamic weights are typically set to 1, and the parameter $\alpha = 0.6$ is used.

\begin{figure}[H]
\centering
\includegraphics[scale=.5]{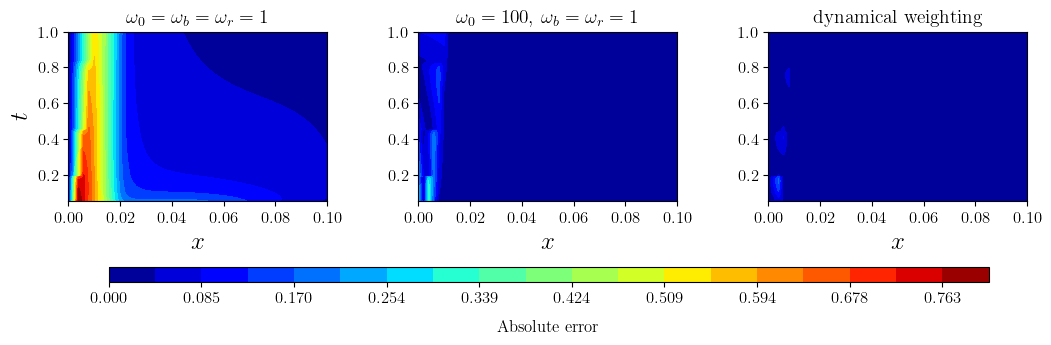}
\caption{Absolute error $|\widehat\theta - \theta|$ for Eq.~\eqref{eq:3.9} at the end of training (i.e. at epoch $=10^5$) for three choices of loss weights.}\label{fig:2}
\end{figure}
\begin{figure}[H]
\centering
\includegraphics[scale=.47]{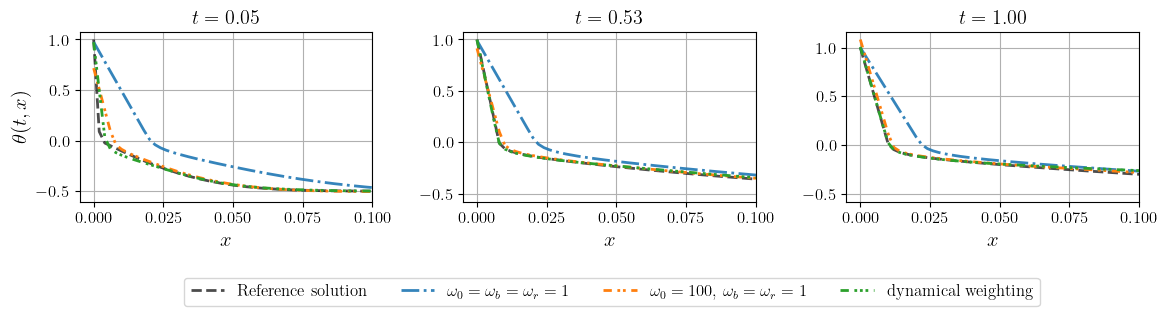}
\caption{Solution of~\eqref{eq:3.9} at the end of training at $t=0.05$, $t=0.53$ and $t=1$.}\label{fig:3}
\end{figure}

From Fig. \ref{fig:2}, we can see that in all considered cases, the predominant portion of the error is located in the region near zero. This is primarily due to the gradual progression of the interface. Indeed, when we set $\Ste=0.005$ and $\Fo=0.01$, we deliberately introduce a movement of the interface with a steep gradient near zero. This type of behavior is difficult to model accurately for a neural network, as they are characterized by a high regularity.\medskip

In Fig. \ref{fig:3}, where the solutions are plotted for $x \in [0, 0.1]$ at various time instants, we can observe distinct behavior for each case. When the weights are set to 1, the model is unable to learn the accurate solution. Furthermore, it seems that the predicted solution get stuck at some intermediate state and cannot be further refined. On the other hand, when $\omega_0=100$ is employed, the initial condition is learned much more effectively, resulting in a significantly more accurate solution at subsequent time instants. In the third scenario, where the learning rate annealing algorithm \cite{wang2021understanding} is applied, we can observe that both the initial and boundary conditions are adhered to, leading to a more precise capture of the interface compared to the other two cases.

\begin{figure}[H]
\centering
\includegraphics[scale=.65]{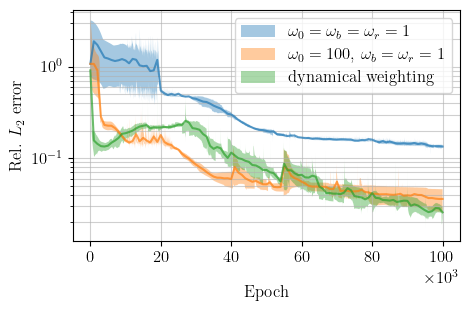}
\caption{Relative $L_2$ error $\|\widehat\theta - \theta\|_2/\|\theta\|_2$ during training for Eq.~\eqref{eq:3.9}.}\label{fig:4}
\end{figure}

Figure \ref{fig:4}, illustrates the significant role of the loss weights in the learning process. For instance, we notice the rapid decrease of the relative $L_2$ error in the first few thousand training epochs for both second (orange) and third (green) cases compared to the first one (blue). After the epoch $60k$, the relative $L_2$ error is barely changed for the first case while it is still considerably decreasing for the other two cases. Furthermore, the error for the third case became smaller than that of the second case. The relative $L_2$ error at the end of training, is $1.341\times 10^{-1} \pm 0.005$ for the first case, $3.59\times 10^{-2}\pm 0.01$ for the second case, and $2.565\times 10^{-2}\pm 0.004$, when the dynamical weighting is performed with a frequency of 1000 iterations.\medskip

\begin{figure}[H]
\centering
\includegraphics[scale=.65]{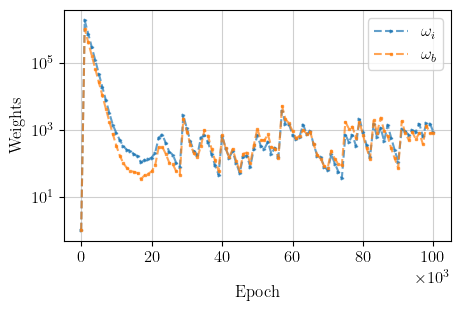}
\caption{Values of the weights during training in the dynamical weighting case for Eq.~\eqref{eq:3.9}.}\label{fig:weights}
\end{figure}

Figure \ref{fig:weights}, shows that the weights are of the same order of magnitude. They reach very large values at the beginning of training to compensate for the dominance of the residual loss function. Over time, these weights gradually decrease and converge to an approximate value of $10^3$.
\pagebreak

\section{Improvements of PINNs for the Stefan problem}

In some cases, the PINNs method can struggle to capture an accurate solution. This difficulty may stem from the choice of various hyperparameters, such as the number of training data, the neural network architecture, the training parameters, etc. Despite careful consideration of these factors, there are instances where the method may still not perform well due to challenges in accurately capturing the underlying physics. In this section, we will explore some strategies aimed at refining the approximation of physics-informed neural networks.

\subsection{Pointwise weighting using soft attention mask functions}

As illustrated in the previous section, the weights within the loss function are a significant factor in the training process. In reality, these weights act as an adjustment to the learning rate during gradient descent for each component of the loss function. In \cite{mcclenny2020self} a pointwise weighting is suggested, to take in account the spatial imbalance of gradients; for every training (or collocation) point $(t, x)$ we assign a weight $\omega(t, x)$. For the problem associated to Eqs. \eqref{eq:3.1}, \eqref{eq:3.2} and \eqref{eq:3.3}, the structure of the loss function becomes:
\begin{equation}
\cal L(\Theta, \omega) = \cal L_0(\Theta, \omega_i) + \cal L_b(\Theta, \omega_b) + \cal L_r(\Theta, \omega_r),
\end{equation}
where $\cal L_r$, $\cal L_0$ and $\cal L_b$ are loss terms corresponding to the PDE term \eqref{eq:3.1}, the initial condition \eqref{eq:3.2}, and the boundary condition \eqref{eq:3.3}, respectively:
\begin{align}
\cal L_r(\Theta) &= \frac{1}{N_r} \sum_{k=1}^{N_r} \,m_r(\omega_r^k)\left|\frac{\partial u_\Theta}{\partial t}(t_r^k, x_r^k) + \cal N[u_\Theta](t_r^k, x_r^k)\right|^2,\\
\cal L_0(\Theta) &= \frac{1}{N_0} \sum_{k=1}^{N_0} \,m_i(\omega_i^k)\left|u_\Theta(0, x_i^k) - g(x_i^k)\right|^2,\\
\cal L_b(\Theta) &= \frac{1}{N_b} \sum_{k=1}^{N_b} \,m_b(\omega_b^k)\left|\cal B[u_\Theta](t_b^k, x_b^k) - h(t_b^k, x_b^k)\right|^2,
\end{align}
where $m_i$, $m_b$ and $m_r$ are nonnegative, differentiable and strictly increasing functions, called mask functions. The weights $\{\omega_i^k\}_k$, $\{\omega_b^k\}_k$ and $\{\omega_r^k\}_k$ are updated by gradient descent side-by-side with the network parameters such that
\begin{equation}
\min_{\Theta}\max_{\omega_i, \omega_b, \omega_r} \cal L(\Theta, \omega_i, \omega_b, \omega_r)   
\end{equation}
is reached.\bigskip

For the Stefan problem \eqref{eq:3.9}, the weighting is specifically applied to the initial and residual points. We employ mask functions $(m_r, m_i, m_b)$ of the kind $w\mapsto \alpha (1 + e^{-\beta (w - m)})^{-1}$, where $\alpha = 1000$, $\beta = 0.1$, and $m = 2$ for the initial condition points, while $\alpha = 1$, $\beta = 1$, and $m = 5$ are chosen for the residual points. The weights are initialized from a uniform distribution in $[0, 1]$ and are updated by gradient descent using the Adam optimizer with a learning rate $\eta_{\rm max} = 10^{-3}$ to maximize the loss function with respect to the weights $w_i$, $w_r$. To minimize the loss function with respect to the model parameters  we use the Adam optimizer for $100k$ iteration with an exponential decay learning rate $t\mapsto \eta_{\rm min}\,\gamma^{t/\kappa}$ where $\eta_{\rm min} = 10^{-3}$, $\gamma = 0.9$, $\kappa = 5000$ and $t$ is the training iteration. All other hyperparameters remain unchanged.\medskip

\begin{figure}[H]
\centering
\includegraphics[scale=.55]{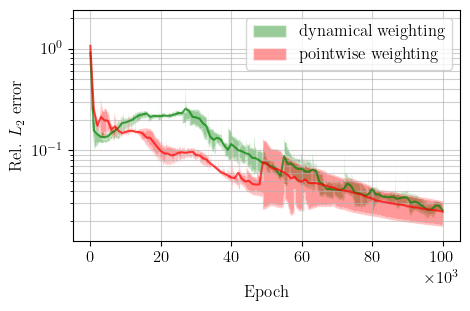}
\hspace{1ex}
\includegraphics[scale=.55]{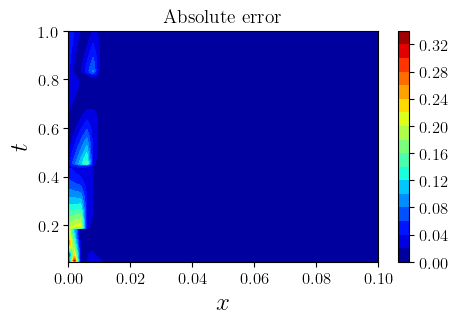}
\caption{Problem \eqref{eq:3.9}. Relative $L_2$ error $\|\widehat\theta - \theta\|_2/\|\theta\|_2$ during training (left). Absolute error $|\widehat \theta - \theta|$ at the end of training (i.e. at epoch $=10^5$) (right).}\label{fig:5}
\end{figure}

\begin{figure}[H]
\centering
\includegraphics[scale=0.47]{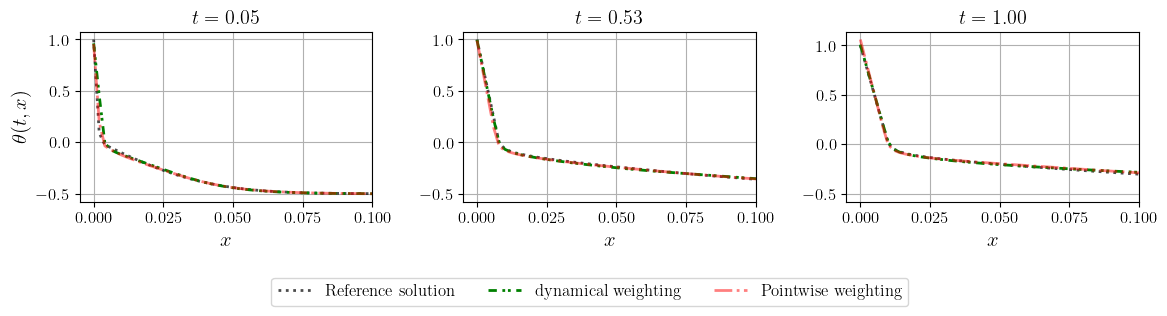}
\caption{Solution of \eqref{eq:3.9} at the end of training for the three cases in $t=0.05$, $t=0.53$ and $t=1$.}\label{fig:6}
\end{figure}

Figure \ref{fig:5}, reveals that during training, the initial $L_2$ error is lower with pointwise weighting than with dynamical weighting. However, beyond epoch $50k$, they converge to similar values (Fig. \ref{fig:6}, shows the solutions at various time instants). The relative $L_2$ error at the end of training, is $2.484\times 10^{-2}\pm 0.007$ for pointwise weighting and $2.565\times 10^{-2}\pm 0.004$ the dynamical weighting.\medskip

\begin{figure}[H]
\centering
\includegraphics[scale=0.6]{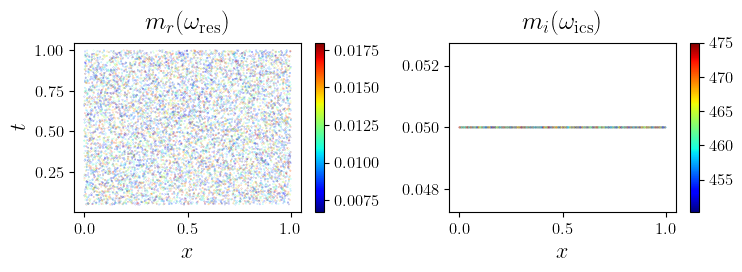}
\includegraphics[scale=0.6]{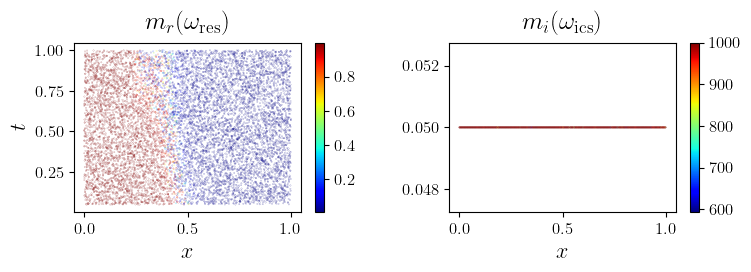}
\caption{Problem \eqref{eq:3.9}. Values taken by the mask function for both the initial and residual points in the beginning of training (top) and at the end of training (bottom).}\label{fig:7}
\end{figure}

Figure \ref{fig:7}, illustrates the multiplicative soft attention mask function values at the beginning and end of training. Initially, when the weights are initialized from a uniform distribution in the range $[0, 1]$, the values corresponding to the residual points are noticeably smaller than those for the initial points. This discrepancy is attributed to the hyperparameters $\alpha$, $\beta$, and $m$, which emphasizes the learning for the initial condition, while constraining the dominance of the residual. At the end of training, we can observe that more attention is needed early in the solution, but not uniformly across the space variable.\medskip

To conclude, this approach can be an alternative to the dynamic weighting method presented in \cite{wang2021understanding}. However, the selection of a mask function for each loss component requires prior knowledge about the challenges involved in fitting the problem.

\subsection{Sequence learning in time}

The conventional PINNs method developed by Raissi et al. \cite{RAISSI2019686} trains the neural network model to predict the solution across the entire space-time domain simultaneously. In certain cases, this can be difficult to learn. An alternative approach proposed in \cite{wight2020solving} involves training the model on sequences of the complete domain. Sequence learning approaches aim to simplify the training optimization problem, and have been proven to be effective \cite{krishnapriyan2021characterizing}.\medskip

In \cite{wight2020solving} Colby et al. considered the sequence $[0, k\Delta t]\times \Omega$, $k=1,\dots, N$ of the time-space domain $[0, T]\times \Omega$, where $N$ is the number of subdomains and $\Delta t = T/N$ the step size. The model is trained firstly on $[0, \Delta t]\times \Omega$, then on $[0, 2\Delta t]\times \Omega$, and so on, gradually increasing the time span up to $[0, N\Delta t]$. To ensure that the solution is well-learned in each time subinterval, one can set a threshold or a maximum number of training iterations. Once the loss function value is below the threshold, or the training has reached the maximum number of iterations, the training procedure begins on the next time interval. In some cases, if the loss function value is still very high after the maximum number of iterations, the time step size $\Delta t$ may need to be reduced.\bigskip

For our computational domain $D = [0.05, 1]\times [0, 1]$, we chose subdomains of the form $D_k = [0.05, 0.05 + k\Delta t]\times [0, 1]$ with uniform time step size $\Delta t = 0.05$. The number of collocation points $N^{(k)}_{r}$ used for $D_k$ is $1000$ if $k = 1$ and $N^{(k - 1)}_{r} + 500$ otherwise. In order to prioritize the initial stages of learning, we determined the number of training iterations $N^{(k)}_{\rm it}$ for each subdomain such that the product $N^{(k)}_{\rm it} \times N^{(k)}_{r}$ is approximately equal to $10^8$.\medskip

\begin{figure}[H]
\centering
\includegraphics[scale=.5]{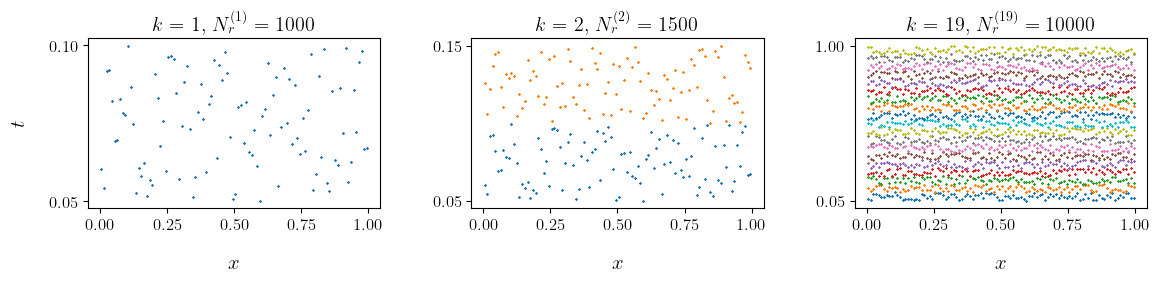}
\caption{Example of collocation points for subdomains $D_k$. Note that the points used in subdomain $D_{k-1}$ are also used in $D_k$.}\label{fig:8}
\end{figure}

Using this set of parameters, the total count of collocation points stands at $10000$, which corresponds to the number used in the previous methods. Furthermore, the maximum number of iterations matches the minimum number of collocation points, ensuring a faster training process and more meaningful comparisons with the previous methods.\medskip

\begin{figure}[H]
\centering
\includegraphics[scale=.6]{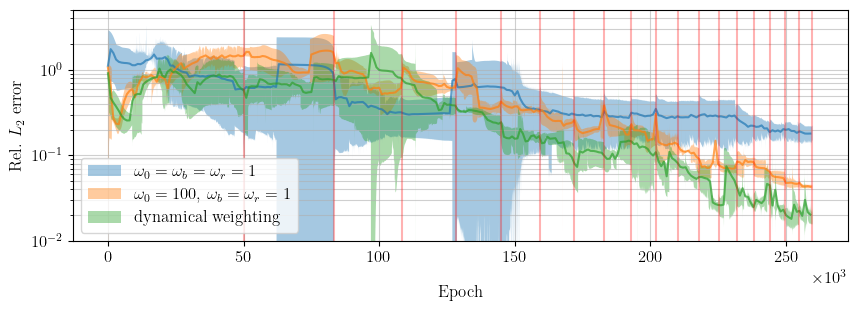}
\caption{Problem \eqref{eq:3.9}. Relative $L_2$ error $\|\widehat\theta - \theta\|_2/\|\theta\|_2$ during training. The red horizontal lines indicate the change of $D_k$ sets.}\label{fig:9}
\end{figure}

Figure \ref{fig:9}, shows the relative $L_2$ error during the training of the model. We notice that in the case when the coefficients of the loss function are set in a static way (i.e. $\omega_0 = \omega_b = \omega_r = 1$ or $\omega_0 = 100$, $\omega_b = \omega_r = 1$), the use of this approach did not improve the results. However in the case where the coefficients are set dynamically the error decreased. The relative $L_2$ error at the end of training, is $1.819\times 10^{-1} \pm 0.036$ when $\omega_0 = \omega_b = \omega_r = 1$, $4.282\times 10^{-2} \pm 0.001$ when $\omega_0 = 100$ and $\omega_b = \omega_r = 1$, and $1.887\times 10^{-2} \pm 0.003$ when the weights are updated dynamically.

\begin{figure}[H]
\centering
\includegraphics[scale=.5]{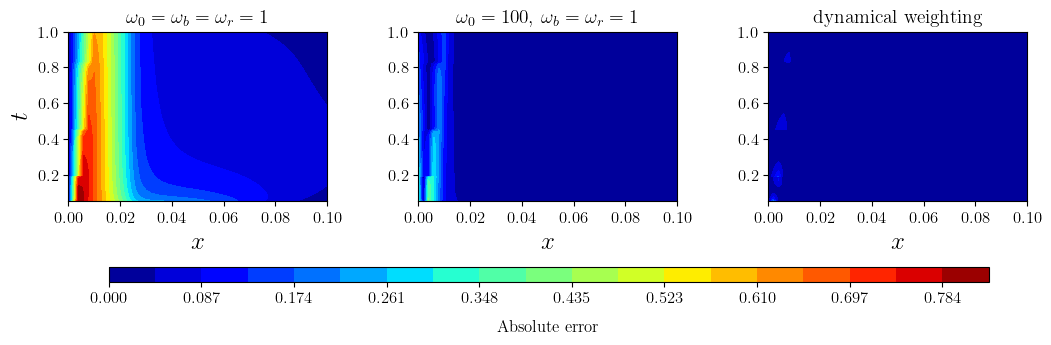}
\caption{Absolute error $|\widehat\theta - \theta|$ for Eq.~\eqref{eq:3.9} when using sequence learning in time.}
\end{figure}

\begin{figure}[H]
\centering
\includegraphics[scale=.5]{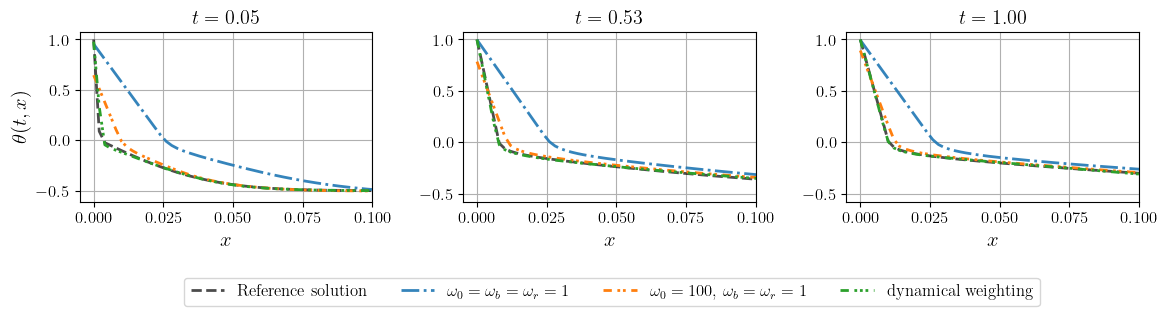}
\caption{Solution of \eqref{eq:3.9} when using sequence learning in time.}
\end{figure}

\section{Conclusion}

We focused on solving Stefan problem using the physics informed neural networks (PINNs) method. We explored various strategies to enhance the approximation capabilities of the PINNs, while maintaining a fixed neural network structure and a consistent number of training and collocation points.\medskip

By choosing a moderate value for the Stefan number (e.g. $\Ste = 0.5$), the Stefan problem can be  solved accurately using the PINNs method, even in the absence of loss function weighting. However, for smaller values (e.g. $\Ste = 0.005$), the solution tends to exhibit less regularity, leading to increased complexity in capturing the interface. In this latter case, the exploration of new strategies becomes indispensable.\medskip

The classical approach employed in \cite{RAISSI2019686} does not account for the imbalance in gradients of the loss component. We have observed that incorporating weights in the loss function proved beneficial and enhanced the accuracy of the approximation. The drawback to manually selecting weights is the ongoing challenge of determining the suitable weights. Considering that, the learning rate annealing algorithm \cite{wang2021understanding} has proven useful in identifying appropriate weights for the loss function. However, based on our experiences, these weights depend on the number of training and collocation points ($N_i$, $N_b$ and $N_r$). 

While weighting can be performed on a pointwise basis, determining the appropriate mask function is not straightforward. The approach of sequence learning in time holds promise for achieving higher accuracy since the optimization problem is confined to a smaller domain. However, a potential drawback, especially with more complex equations, is that if the solution deviates within a time interval, subsequent intervals will propagate this error. It is important to learn the solution well on a time interval before progressing to the next one. This method may also be helpful when working on problems with larger time domains.\medskip

In general, a single approach can yield different results for various problems. In this work, we aimed to explore different methods to solve the Stefan problem. We did not examine the impact of network size, learning rate, or the number of training data on the accuracy of the approximation. Considering these parameters could lead to more precise solutions, and this will be one of the directions for our future research.\smallskip

The Stefan problem studied in this paper is a classical phase change problem where the motion of the liquid phase, described by the Navier-Stokes equations and natural convection is ignored. In future studies, we will focus on a more realistic model of phase change problems, specifically the Navier-Stokes-Boussinesq equations.

\section*{Acknowledgement}

This work has been supported by the regional project \href{https://www.datalab-normandie.fr/}{DATALAB}.
The authors also acknowledge financial support from the French ANR grant  {ANR-23-CE23-0020} (project PINNterfaces). Part of this work used computational resources provided by CRIANN (Centre R{\'e}gional Informatique et d'Applications Num{\'e}riques de Normandie).

\bibliographystyle{unsrt}
\bibliography{stefan}

\end{document}